\magnification=1200

\input xy
\xyoption {all}

\centerline {\bf A SERIES OF SMOOTH IRREGULAR VARIETIES}
\vskip 0.5truecm
\centerline {\bf IN PROJECTIVE SPACE}
\vskip 0.5truecm

\centerline {Ciro Ciliberto and Klaus Hulek}

\vskip 2.0 truecm

\noindent {\bf Abstract. }\par
One of the simplest examples of a smooth, non degenerate surface
in ${\bf P}^ 4$
is the quintic elliptic scroll. It can be constructed from an elliptic
normal curve $E$ by joining every point on $E$ with the
translation of this point by a non-zero $2$-torsion point.
The same construction can be
applied when $E$ is replaced by a (linearly normally embedded)
abelian variety $A$. In this paper we ask the question when the resulting
scroll $Y$ is smooth. If $A$ is an abelian surface embedded by
a line bundle $L$ of type $(d_1,d_2)$ and $r=d_1d_2$, then we prove
that for general
$A$ the scroll $Y$ is smooth if $r$ is at least $7$ with the one exception
where $r=8$ and the $2$-torsion point is in the kernel $K(L)$ of $L$.
In this case $Y$
is singular. The case $r=7$ is particularly interesting, since then
$Y$ is a smooth threefold in ${\bf P}^ 6$  with irregularity $2$.
The existence of this
variety seems not to have been noticed before.
One can also show that the case of the quintic elliptic scroll
and the above case are the only
possibilities where $Y$ is smooth and the codimension
of $Y$ is at most half the dimension of the surrounding projective space. \par

\vskip 0.2truecm

 \par

\vskip 0.4truecm

\noindent {\bf 0. Introduction. }\par

\vskip 0.2truecm

One of the simplest examples of a smooth, non degenerate surface in ${\bf
P}^ 4$ is the
quintic elliptic scroll $Y$. Its construction goes as follows. Let $A$ be
an elliptic normal
curve of degree $5$ in ${\bf P}^ 4$ and let $\epsilon$ be a non zero point
of order two on
$A$. Then the union of all the lines joining pairs of points of type
$x$ and $x+\epsilon$
on
$A$ is an elliptic quintic scroll.\par

Exactly the same construction can be repeated starting from any abelian
variety $A$ of
dimension $n$, with $A$ linearly normally embedded in a projective space
${\bf P}^ N$ via
a very ample line bundle ${\cal L}$, and from any non trivial point
$\epsilon\in A$ of order
two. We investigate this construction in the present paper. In this way we
get a
scroll
$Y$ of dimension $n+1$ in ${\bf P}^ N$ related to the above data
$(A,\epsilon,{\cal L})$ and
the first interesting question is: when is
$Y$ smooth? It is well known that this is the case if $n=1$ and $N\geq 4$.
So the
next interesting case is that of surfaces, i.e. $n=2$, embedded in ${\bf
P}^ {r-1}$ via
a $(d_1,d_2)$-polarization, with $r=d_1\cdot d_2$. If $r\leq 6$ there is no
hope for $Y$ to
be smooth because of Lefschetz's hyperplane section theorem. So the
question becomes relevant
as soon as $r\geq 7$. In fact the main part of this paper is devoted to
proving that if $A$
is general in its moduli space (it is enough to assume that
$End(A)\simeq {\bf Z}$ or $NS(A)\simeq {\bf Z}$
depending on the case under consideration), and if
$r\geq 7$ and $r\not=8$, then $Y$ is smooth. This is particularly remarkable
in the case
$r=7$, since
$Y$ is then an irregular, codimension three manifold in ${\bf P}^ 6$, whose
existence
does not seem to have been previously noticed. As we remark at the end of
\S 2, for no
other dimension of
$A$, but
$1$ and $2$, and
$N=4$ and $N=6$ respectively, $Y$ can be smooth of codimension $c\leq
{N\over 2}$
in
${\bf P}^ N$. The case $d=8$ is also interesting. If $A$ is a general
abelian surface
with a polarization of type $(1,8)$, then $Y$ is smooth, unless the
translation by
the point $\epsilon$ of order two fixes the polarization, in which case $Y$
is singular.
If the polarization is of type $(2,4)$, then the translation by $\epsilon$
automatically
fixes the
polarization and $Y$ is again singular.\par

The paper is organised as follows. In \S 1 we present the construction of a
suitable
projective bundle $X$ over ${\overline A}=A/\epsilon$ which maps to $Y$ via its
tautological line bundle. In \S 2 we prove that this map is finite and we
compute the
double point cycle of the composite map of $X\to Y$ with a general
projection in a ${\bf P}^
l$, with $n+1\leq l\leq N$. From \S 3 on we restrict our attention to the
case of abelian
surfaces. In particular in
\S 3 we prove that $Y$ is smooth as soon as $r\geq 10$. This comes as a
consequence of the
fact that, in this situation, if $A$ is general enough, then it has no
quadrisecant plane. A property which, in turn, follows as an application of
Reider's method. Finally in
\S 4 we prove that $Y$ is smooth if
$r=7$, in \S 5 we analyse the case
$r=8$ and in \S 6 the case $r=9$.
The idea for the proof that $Y$ is smooth and the tools we use in
the cases
$r=7$,
$r=8$ and the polarization is of type $(1,8)$ with $\epsilon$ not fixing
it, and $r=9$ and the
polarization is the triple of a principal polarization (which is the only
critical case for
$r=9$) are the same: we first bound
dimension and degree of the possible singular locus of
$Y$ by using geometric arguments and the double point formula, then we use
the action of the
Heisenberg group
to give a lower bound for the degree of the singular locus, finally
contradicting the
previous estimate.\par

\vskip 1.0 truecm

\noindent {\bf Acknowledgments.} The present collaboration took place in the
framework of
the HCM contract AGE (Algebraic Geometry in Europe), no. ERBCHRXCT940557.\par
We are very grareful to the referee who has not only made a number of
suggestions which have improved the presentation of the paper, but who has
also pointed out some inaccuracies such as that the formula given in
proposition (2.3) was incorrect as it was stated in the first version of
this paper. The referee's comments also lead to a shorter
proof of theorem (5.1)

\vskip 1.0 truecm

\noindent {\bf 1. Some projective bundles over abelian varieties.}\par

\vskip 0.2truecm

Let $A$ be an abelian variety of dimension $n$ with a polarization
$\Theta\in NS(A)$ of type
$(d_1,...,d_n)$ with $d_1|...|d_n$ (our general reference for the theory of
abelian
varieties will be [LB]). Let us take a non trivial point
$\epsilon\in A$ of order two.\par

Let $K(\Theta)$ be the kernel of the isogeny $\lambda_\Theta: A\to {\hat
A}=Pic^0(A)$
determined by the polarization. Recall that $K(\Theta)\simeq ({\bf
Z}_{d_1}\times ...\times
{\bf Z}_{d_n})^ 2$ and that, if $d_1$ is even, then $\Theta$ is divisible
by two in
$NS(A)$ and every point of order two of
$A$ is an element of $K(\Theta)$. \par

Let ${\cal L}$ be a line bundle on $A$ representing $\Theta$. Then we
have:\par

$$t_\epsilon^*{\cal L}\simeq {\cal L}\otimes {\cal L}_0$$

\noindent where $t_x$ is the translation by a point $x\in A$ and
${\cal L}_0\in Pic^0(A)$ is the
point of order two given by $\lambda_\Theta(\epsilon)$. Hence ${\cal L}_0$
is trivial if
and only if $\epsilon\in K(\Theta)$.\par

Let $\overline A$ be the quotient $A/\epsilon$ and let $\pi: A\to \overline
A$ be the quotient
map, which is an isogeny of degree $2$. If $\epsilon\in K(\Theta)$, then
there is a line
bundle $\overline {\cal L}$ on $\overline A$, such that $\pi^*(\overline
{\cal L})={\cal
L}$. The line bundle $\overline {\cal L}$ represents a polarization
$\overline \Theta$ on
$\overline A$, of type $(\overline {d_1},...,\overline {d_n})$, such that:\par

$$2\cdot \overline {d_1}\cdot ...\cdot \overline {d_n}={d_1}\cdot ...\cdot
{d_n},$$

\noindent a relation which is obtained from $\Theta^n=\pi^*({\overline
\Theta})^n=2{\overline
\Theta}^n$. In particular, if
${d_1}= ...= {d_{n-1}}=1$,
$d_n=d$, then
$d$ is even and $\overline \Theta$ is of type $(1,...,1,{d\over 2})$.

One has:\par

$$\pi_*{\cal O}_A={\cal O}_{\overline A}\oplus {\cal M}_1$$

\noindent where ${\cal M}_1$ is a non trivial $2$-torsion point in
$Pic^0(\overline A)$. The
induced map $\pi^*: Pic^0(\overline A)\to Pic^0(A)$ is also an isogeny of
degree $2$, whose
kernel is generated by ${\cal M}_1$. Therefore we have two line bundles
${\cal M}_2, {\cal
M}_3\in Pic^0(\overline A)$ such that:\par

$$\pi^* ({\cal M}_2)=\pi^*({\cal M}_3)={\cal L}_0$$

\noindent and one has ${\cal M}_2={\cal M}_3\otimes {\cal M}_1$. The
elements ${\cal M}_i$,
$i=1,2,3$, and the trivial bundle form a subgroup $G$ of $Pic^0(\overline
A)$, which is the
inverse
image via $\pi^*$ of the subgroup generated by ${\cal L}_0$.

 We have the:\par

\vskip0.4truecm\noindent {\bf
Lemma} 1.1 {\sl If $\epsilon\in K(\Theta)$ then $G$ is the group of order
two generated
by ${\cal M}_1$. Otherwise
$G$ is isomorphic to
${\bf Z}_2\times {\bf Z}_2$}.

\vskip0.2truecm\noindent {\it Proof.} The first assertion is clear.
To prove the second assertion let
$L$ be a lattice which defines $A = {\bf C}^g/L$.
Then the point $\epsilon$ is represented by an element $e \in {1\over 2}L$.
The fact that
$\epsilon$ is not in $K(\Theta)$ is equivalent to the existence of some
element $f \in L$
such that for the pairing defined by the polarisation
$(e,f)= {1\over 2}$ mod ${\bf Z}$.
The lattice
$\overline L$ which defines the quotient $\overline A$ is the lattice
generated by $L$ and
$e$. We denote by $L^{\vee}$ the dual lattice of $L$. This defines the dual
variety $Pic^0(A)={\bf C}^g/L^{\vee}$ of $A$.
The element $e$ is not contained in $L^{\vee}$ and
represents the line
bundle ${\cal L}_0$ in $Pic^0(A)$.
Similarly $f \notin {\overline
L}^{\vee}$, but
$2f \in {\overline L}^{\vee}$. The element $f$ corresponds to the line bundle
${\cal M}_1$ in $Pic^0(\overline A)$ which is $2$-torsion.
The element $e$ also defines a line
bundle in
$Pic^0(\overline A)$, whose pullback to $A$ is ${\cal L}_0$ and this
corresponds to
${\cal M}_2$ or ${\cal M}_3$ .
The claim follows if we can show that ${\cal M}_2$ or ${\cal M}_3$ is
$2$-torsion.
But this follows since $2e \in {\overline L}^{\vee}$.
(A different proof will follow from proposition
(1.5) below (see remark (1.6))).
$\diamondsuit$\par

\vskip0.4truecm

Let us set:

$${\cal E}= {\cal E}(A,\epsilon,{\cal L})=\pi_*{\cal L}.$$

\noindent This is a rank $2$ vector bundle on ${\overline A}$, and we can
consider the associated
projective bundle:\par

$$X=X(A,\epsilon,{\cal L})={\bf P}({\cal E})={\bf Proj}(\oplus_{i=0}^\infty
Sym^i({\cal
E}))$$

\noindent with its tautological line bundle ${\cal O}_X(1)$ and its
structure map $p: X\to
{\overline A}$.  We will denote by $F$ a fibre of $p$ and by $H$ a divisor
in $|{\cal
O}_X(1)|$. We will use the same notation to denote their classes in the
homology ring
of $X$.\par

If $\epsilon\in K(\Theta)$ then, since $\pi^*(\overline {\cal L})={\cal
L}$, the projection formula tells us that:\par

$${\cal E}=\overline {\cal L}\oplus (\overline {\cal L}\otimes {\cal M}_1).$$

\noindent By contrast, as we shall see in a moment, if $\epsilon\not\in
K(\Theta)$, then the
bundle ${\cal E}$ in general does not split.\par

The natural map $\pi^{*} \cal E \to \cal L$ defines an inclusion $i$ of $A$
into $X$ such that $H$ restricts to $\cal L$ on $i(A)$.
The image $i(A)$
is a $2$-section over $\overline A$. If there is no danger of confusion we
shall
denote $i(A)$ also by $A$. Let
$$\tilde X={\bf P}(\pi^*{\cal E}).$$
We have a natural \'etale map $f: \tilde X\to X$ of degree $2$.
The inverse image of the $2$-section $A$ in $X$
under the map $f$ consists of $2$ sections
of $\tilde X$ corresponding to the $2$ projections $\pi^{*} \cal E \to \cal L$
and $\pi^{*} {\cal E} \to t^*_{\epsilon} {\cal L}$ whose existence follows
from the construction of $\cal E$.
This
shows that $\pi^{*} \cal E$ splits, more precisely

$$ \pi^*{\cal E}\simeq {\cal L}\oplus t^*_\epsilon{\cal L}= {\cal L}\oplus
({\cal
L}\otimes {\cal L}_0). \eqno (1) $$
Notice
that, if $\epsilon\in K(\Theta)$, then $\tilde X=A\times {\bf P}^1$ is
trivial.\par

\vskip0.4truecm\noindent {\bf
Lemma} 1.2 {\sl One has:\par
\medskip

\noindent (i)   ${\cal O}_X(1)_{|A}\simeq {\cal L}$;\par

\noindent (ii) ${\cal O}_A(A)\simeq {\cal L}_0$;\par

\noindent (iii) ${\cal O}_X(-2A)\simeq {\cal O}_X(2K_X)$.\par
\medskip

\noindent Moreover, if $\epsilon\not\in K(\Theta)$ then:\par
\medskip

\noindent (iv) there is no section $A'$ of $X$ over ${\overline A}$ which
is disjoint from
$A$.\par
\medskip

\noindent Hence if $\epsilon\not\in K(\Theta)$ and $A$ does not contain
elliptic curves, then
${\cal E}$ does not split.}

\vskip0.2truecm\noindent {\it Proof.} (i) follows by the definition of the
tautological
bundle.\par

(ii) There are two
sections $\tilde A_1$ and $\tilde A_2$ of $\tilde X$ over $A$,
which map both isomorphically to $A$
via $f$. These sections correspond to the splitting of $\pi^*({\cal E})={\cal
L}\oplus ({\cal L}\otimes {\cal L}_0)$. Since the normal bundle of both
these sections in
$\tilde X$ is given by ${\cal L}_0$, and since $f$ is \'etale, we have the
assertion.\par

(iii) Since $A\cdot F=-K_X\cdot F=2$, there is a line
bundle ${\cal M}$ on ${\overline A}$ such that ${\cal O}_X(-A)\simeq {\cal
O}_X(K_X)\otimes
\pi^*({\cal M})$. On the other hand, by adjunction, one has ${\cal
O}_A(-A)\simeq {\cal
O}_A(K_X)$. This implies that either
${\cal M}\simeq {\cal O}_{\overline A}$ or ${\cal M}\simeq
{\cal M}_1$. This immediately yields the assertion. \par

(iv) Suppose $A'$ is disjoint from $A$. Then $A'$ would pull back to a
section $\tilde {A'}$ of $\tilde X$, disjoint from $\tilde A_1$
and $\tilde A_2$, which would give
another way of splitting
$\pi^*({\cal E}\otimes {\cal L}^*)$. This is impossible under the assumption
$\epsilon\not\in K(\Theta)$.
\par

If ${\cal E}$ splits, we have two sections $A'$ and $A''$ of $X$ which do
not meet. However
they both meet $A$ and they must cut out two divisors $C'$ and $C''$
on $A$ which do not meet each
other. Hence $C'$ and $C''$ are pull-backs from an elliptic curve
and so $A$ contains an elliptic curve.
$\diamondsuit$\par

\vskip0.4truecm Notice that the map $p:X \to {\overline A}$ induces an
isomorphism
$p^*:Pic^0({\overline A})\to Pic^0(X)$. We will identify $Pic^0({\overline
A})$ and
$ Pic^0(X)$ using $p^*$. We
have the:\par

\vskip0.4truecm\noindent {\bf
Proposition} 1.3 {\sl Let $\eta$ be an element of $Pic^0(X)$. One has:\par

\noindent (i) if $\epsilon\in K(\Theta)$ then $h^0(X,{\cal
O}_X(A)\otimes
\eta)=0$ unless $\eta={\cal O}_X$, in which case $h^0(X,{\cal O}_X(A))=2$,
and $\eta = {\cal
M}_1$, in which case $h^0(X,{\cal O}_X(A)\otimes {\cal M}_1)=1$;\par

\noindent (ii) if $\epsilon\not\in K(\Theta)$ then $h^0(X,{\cal O}_X(A)\otimes
\eta)=0$ unless $\eta={\cal O}_X, {\cal M}_2, {\cal M}_3$,
in which cases $h^0(X,{\cal O}_X(A)\otimes
\eta)=1$.}

\vskip0.2truecm\noindent {\it Proof.} We have:\par

$$p_*{\cal O}_X(1)\simeq {\cal E},\quad  p_*{\cal O}_X(2)\simeq Sym^2{\cal
E}.$$

Therefore, by using $(1)$ we have:\par

$$ \pi^*p_*{\cal O}_X(2)\simeq Sym^2\pi^*{\cal E}\simeq {\cal L}^{\otimes
2}\oplus
{\cal L}^{\otimes 2}\oplus ({\cal L}^{\otimes 2}\otimes {\cal L}_0). \eqno
(2)$$

\noindent Moreover, since $A\cdot F=2$, there is a
line bundle ${\cal N}$ on ${\overline A}$ such
that:\par

$$ {\cal O}_X(A)\simeq {\cal O}_X(2)\otimes p^*{\cal N}. \eqno (3)$$

\noindent Hence, by the projection formula, one has:\par

$$ \pi^* p_*({\cal O}_X(A)\otimes \eta)\simeq \pi^*(p_*{\cal O}_X(2)\otimes
({\cal N}\otimes
\eta))\simeq \pi^* p_*{\cal O}_X(2) \otimes \pi^* ({\cal N}\otimes
\eta). \eqno (4)$$

\noindent From lemma (1.2, i, ii) and restricting (3) to $A$, we obtain:\par

$$\pi^*{\cal N}^*\simeq {\cal L}^{\otimes 2}\otimes {\cal L}_0. \eqno (5) $$

\noindent Now, by (2), (4) and (5), we get:\par

$$\pi^*p_*({\cal O}_X(A)\otimes \eta)
\simeq ({\cal L}_0\oplus {\cal L}_0\oplus {\cal O}_A)\otimes
\pi^*\eta$$

\noindent and therefore if $\epsilon\in K(\Theta)$:\par

$$ h^0(A,\pi^*p_*{\cal O}_X(A\otimes \eta)) =
\cases {
  0 & if $ \pi^*\eta \neq {\cal O}_A $ \cr
  3 & if $ \pi^*\eta = {\cal O}_A $, i.e.\
         $ \eta\simeq {\cal O}_{\overline A},{\cal M}_1 $ \cr
}
\eqno (6) $$

\noindent whereas:

$$ h^0(A,\pi^*p_*{\cal O}_X(A\otimes \eta)) =
\cases {
  0 & if $ \pi^*\eta \neq {\cal O}_A,{\cal L}_0 $ \cr
  1 & if $ \pi^*\eta = {\cal O}_A $, i.e.\
         $ \eta\simeq {\cal O}_{\overline A}, {\cal M}_1 $ \cr
  2 & if $ \pi^*\eta = {\cal L}_0 $, i.e.\
         $ \eta \simeq {\cal M}_2, {\cal M}_3 $ \cr
}
\eqno (6') $$

\noindent otherwise. Notice that:\par

$$h^0(A,\pi^*p_*{\cal O}_X(A\otimes \eta))=h^0({\overline A},
\pi_*\pi^*p_*{\cal
O}_X(A\otimes \eta))=$$
$$=h^0({\overline A},p_*{\cal O}_X(A\otimes \eta))+h^0({\overline
A},(p_*{\cal O}_X(A\otimes
\eta))\otimes {\cal M}_1)=$$
$$= h^0(X,{\cal O}_X(A\otimes \eta))+h^0(X,{\cal O}_X(A\otimes
\eta)\otimes {\cal M}_1).$$

\noindent Then, by $(6)$, resp. $(6')$
if $\eta\neq {\cal O}_{\overline A}, {\cal M}_i, i=1,2,3$, we
have: \par

$$h^0(X,{\cal O}_X(A\otimes \eta))+h^0(X,{\cal O}_X(A\otimes
\eta)\otimes {\cal M}_1)=0. $$

\noindent in particular $h^0(X,{\cal O}_X(A\otimes \eta))=0$.\par

Let $\epsilon\in
K(\Theta)$. If $\eta = {\cal O}_{\overline A}, {\cal M}_1$, we find:\par

$$h^0(X,{\cal O}_X(A))+h^0(X,{\cal O}_X(A)\otimes {\cal M}_1)=3. $$

\noindent We claim that $h^0(X,{\cal O}_X(A))=2$
and $h^0(X,{\cal O}_X(A)\otimes {\cal M}_1)=1$. Let $\tilde A$
be a trivial section of $\tilde X=A\times {\bf P}^1$ over $A$. Of course
$|\tilde A|$ is
a base point free pencil on $\tilde X$. The image of this pencil under the
map $f: \tilde X\to
X$ is a system of divisors on $X$ which is contained in a linear system.
Since $A$ is in this
system, we see that $h^0(X,{\cal O}_X(A))\geq 2$. on the other hand we
cannot have
$h^0(X,{\cal O}_X(A))\geq 3$ because of lemma (1.2, ii). Hence the
assertion follows,
proving (i).\par

Let now $\epsilon\not\in
K(\Theta)$. If $\eta = {\cal O}_{\overline A},
{\cal M}_1$, we find:\par

$$h^0(X,{\cal O}_X(A))+h^0(X,{\cal O}_X(A)\otimes {\cal M}_1)=1, $$

\noindent and since $h^0(X,{\cal O}_X(A))\geq 1$, we have $h^0(X,{\cal
O}_X(A))=1$ and
$h^0(X,{\cal O}_X(A)\otimes {\cal M}_1)=0$. Finally, if $\eta = {\cal M}_i,
i=2,3$, we
have:\par

$$h^0(X,{\cal O}_X(A)\otimes {\cal M}_2)+h^0(X,{\cal O}_X(A)\otimes {\cal
M}_3)=2. $$

We claim that both summands are smaller than $2$. Otherwise the linear
system, say, $|{\cal
O}_X(A)\otimes {\cal M}_2|$ would be a pencil, and therefore we would find
an element of
it meeting $A$. But by lemma (1.2, ii), the restriction of ${\cal
O}_X(A)\otimes {\cal M}_2$ to $A$ is trivial. This would yield that
$A$ itself is an element of the pencil, implying $h^0(X,{\cal M}_2) > 0$
and hence that
${\cal M}_2$ is trivial on
$X$, a contradiction. In
conclusion we have:\par

$$h^0(X,{\cal O}_X(A)\otimes {\cal M}_2)=h^0(X,{\cal O}_X(A)\otimes {\cal
M}_3)=1 $$

\noindent which finishes our proof. $\diamondsuit$\par

\vskip0.4truecm\noindent {\bf
Remarks} 1.4  (i)  First we consider the case $\epsilon\in K(\Theta)$. We
reconsider the
relation between the pencils $|\tilde A|$ on $\tilde X$ and $|A|$ on $X$.
The map $f$ sends
each element of $|\tilde A|$ to an element of $|A|$. As we saw in the proof
of lemma (1.2,
ii), we have $f^*(A)=\tilde A_1+\tilde A_2$. Hence $f$ is two-to-one
between $|\tilde A|$
and $|A|$. This means that all, but two, elements of $|A|$ are smooth,
irreducible,
isomorphic to $A$, and that there are two elements of $|A|$ of type
$2{\overline A}^+$,
$2{\overline A}^-$
with ${\overline A}^{\pm}$ sections of $X$ over ${\overline A}$. One moment of
reflection shows that these two sections, which do not meet, correspond to
the splitting of
${\cal E}$. Of course
${\overline A}^{\pm}$ are isomorphic to
${\overline A}$ and one has ${\cal O}_{{\overline A}^{\pm}}({\overline
A}^{\pm})\simeq {\cal
M}_1$. In addition $2{\overline A}^{\pm}\equiv A$ but of course ${\overline
A}^+\not\equiv {\overline A}^-$. Hence ${\overline A}^+-{\overline A}^-$
gives a
non trivial point of order two in $Pic^0(X)\simeq Pic^0(\overline A)$. By
restricting to
${\overline A}^{\pm}$, we see that this point of order two is ${\cal M}_1$.
Hence ${\cal O}_X({\overline A}^+)\simeq  {\cal O}_X({\overline
A}^-)\otimes {\cal M}_1$ and
therefore  ${\cal O}_X(A)\simeq {\cal O}_X(2{\overline A}^+)\simeq
{\cal O}_X({\overline A}^++{\overline A}^-)\otimes {\cal M}_1$, whence
${\cal O}_X(A)\otimes {\cal M}_1\simeq
{\cal O}_X({\overline A}^++{\overline A}^-)$, which fully explains the
meaning of part
(i) of proposition (1.3).\par
Notice that all the smooth abelian varieties in $|A|$ play a symmetric role
in the
construction of $X$ and of its tautological line bundle.\par
One more obvious remark. Let ${\overline A}^-$ correspond to the quotient
${\cal E}\to
{\overline {\cal L}}$ and ${\overline A}^+$ to the quotient ${\cal E}\to
{\overline {\cal L}\otimes {\cal M}_1}$. Then ${\cal O}_{{\overline
A}^-}(1)\simeq
{\overline {\cal L}}$ and ${\cal O}_{{\overline A}^+}(1)\simeq
{\overline {\cal L}\otimes {\cal M}_1}$.\par
(ii) Now we take up the case $\epsilon \not\in K(\Theta)$.
Consider the
varieties $A_2$, $A_3$ which are the unique divisors in the linear systems
$|{\cal O}_X(A)\otimes {\cal M}_3|$, $|{\cal O}_X(A)\otimes {\cal M}_2|$,
respectively. As we
saw in the proof of proposition (1.3), we have
$A\cap A_2=A\cap A_3=\emptyset$. Then, by lemma (1.2, iv), $A_2$ and $A_3$
are irreducible. We shall see in proposition (1.5) that these varieties are
smooth abelian. We also set $A_1=A$.

\vskip0.4truecm\noindent {\bf
Proposition} 1.5 {\sl One has:\par

\noindent (i) if $\epsilon\in K(\Theta)$, then $h^0(X,{\cal O}_X(2A))=3$;\par

\noindent (ii) if $\epsilon\not\in K(\Theta)$, then
$h^0(X,{\cal O}_X(2A))=2$. Moreover the pencil $|2A|$ has exactly $3$ singular
elements namely $2A_i$ for $i=1,2,3$. All other elements $D$ in $|2A|$ are
smooth abelian. The reduced varieties $A_i$ are also smooth abelian.  }

\vskip0.2truecm\noindent {\it Proof.} In case (i) the linear system is
composed with
the pencil $|A|$, hence the assertion. \par

Let us consider case (ii). Since
$2A\equiv 2A_2\equiv 2A_3$, it is clear that
$h^0(X,{\cal O}_X(2A))\geq 2$. Suppose $h^0(X,{\cal O}_X(2A)) = r+1\geq 3$.
Then the linear system
$|2A|$ would have dimension $r\geq 2$. Moreover ${\cal O}_A(2A)$ is
trivial. Therefore the
linear system $|A|$ would have dimension at least $r-1\geq 1$,
contradicting proposition (1.3). We now look at the pencil $|2A|$.
We have already seen that the $A_i$ are irreducible. The same is true for
the elements $D$.
In fact $D$ does not meet $A$ and
by lemma (1.2, iv) $D$ cannot have a component which is a section.
In addition, there cannot be a non trivial component which is a $2$-section
either,
because such a
$2$-section would be numerically equivalent to $A$, hence would be equal to
$A_2$ or
$A_3$, which is not possible since $A_2 + A_3$ is not equivalent to $2A$.
It follows from adjunction that the square
of the dualizing sheaf $\omega_{A_i}$ is trivial and that $\omega_D$ is
trivial. Our assertion follows if we can show that the projection onto
$\overline A$ defines an \'etale $2:1$ cover from $A_i$ to $\overline A$,
resp. an \'etale $4:1$ cover from $D$ to $\overline A$. To see this we look
at the pencil of degree $4$ cut out by $|2A|$ on each ruling. This is
base point free and has at least $3$ singular elements consisting of
$2$ double points each, corresponding to the $2$-sections $A_i$. By the
Hurwitz formula there can be no worse singularities and this gives the claim.
$\diamondsuit$\par

\vskip0.4truecm\noindent {\bf
Remarks} 1.6 (i) Assume that $\epsilon\not\in K(\Theta)$. Then
we have just seen that $A_2$ and $A_3$ are smooth abelian
varieties isogenous to ${\overline A}$ via
the degree $2$ maps $\pi_2$, $\pi_3$ induced by $p$. In addition
we have $A_2\cap A_3=\emptyset$.
Hence $\pi_2^*{\cal M}_2\simeq {\cal O}_{A_2}$ and $\pi_3^*{\cal M}_3\simeq
{\cal O}_{A_3}$,
which gives another proof of lemma (1.1) in the present case.\par

(ii) In this situation the projection $p$ induces on every
smooth element $D\in |2A|$  an isogeny
$\delta: D\to {\overline A}$ of degree $4$. We have just seen that
$\pi_2^*{\cal M}_2\simeq {\cal O}_{A_2}$ and $\pi_3^*{\cal M}_3\simeq
{\cal O}_{A_3}$. It follows from this that
$\delta^* ({\cal
M}_i)$ is trivial for every $i=1,2,3$. Hence $D$ is constant in
moduli and it is the
unique degree $4$ cover of ${\overline A}$ with this property.
We also remark that, in view of the above description, the isogeny
$\delta$ factors through
degree $2$ isogenies
$\delta_i: D\to A_i$,
$i=1,2,3$.
\par

(iii) We consider the line bundles
${\cal L}_i:={\cal
O}_X(1)_{|A_i}$, and the corresponding polarizations $\Theta_i$,
$i=1,2,3$. Of course ${\cal E}\simeq {\pi_i}_*{\cal L}_i$, $i=1,2,3$, and
the abelian varieties
$A_i$ play a symmetric role in the construction of $X$ and of its
tautological line bundle
${\cal O}_X(1)$.

\par

\vskip0.4truecm\vskip0.4truecm\noindent {\bf
Proposition} 1.7 {\sl One has:\par

\noindent (i) if $\epsilon\in K(\Theta)$, then ${\cal O}_X(K_X)\simeq {\cal
O}_X(-A)$;\par

\noindent (ii) if $\epsilon\not\in K(\Theta)$, then
${\cal O}_X(K_X)\simeq {\cal
O}_X(-A_i)\otimes {\cal M}_i$, for $i=1,2,3$. }

\vskip0.2truecm\noindent {\it Proof.} (i) By the proof of lemma (1.2), we
know that either
${\cal O}_X(K_X)\simeq {\cal
O}_X(-A)\otimes {\cal M}_1$ or ${\cal
O}_X(K_X)\simeq {\cal
O}_X(-A)$. The assertion follows by restricting to ${\overline A}^{\pm}$.\par

(ii) As above the proof of lemma (1.2) tells us that either
${\cal O}_X(K_X)\simeq {\cal
O}_X(-A_i)\otimes {\cal M}_i$ or ${\cal
O}_X(K_X)\simeq {\cal
O}_X(-A_i)$. Suppose that ${\cal
O}_X(K_X)\simeq {\cal
O}_X(-A)$.
Then the adjunction formula tells us
that ${\cal
O}_{A_2}\simeq {\cal
O}_X(K_X+A_2)\otimes {\cal
O}_{A_2}\simeq {\cal
O}_X(A_2-A_1)\otimes {\cal
O}_{A_2}\simeq {\cal M}_3\otimes {\cal
O}_{A_2}$, which is a contradiction, since only ${\cal M}_2\otimes {\cal
O}_{A_2}$ is trivial on $A_2$. $\diamondsuit$

\vskip0.4truecm

Let us now consider the action of $K(\Theta)$ on $A$. We will assume
$\Theta$ is a
primitive polarization, i.e. it is an indivisible element of $NS(A)$, of type
$(d_1,d_2,...,d_n)$. This is equivalent to $d_1=1$.

\vskip0.4truecm\vskip0.4truecm\noindent {\bf
Lemma} 1.8 {\sl Let us suppose that the Neron-Severi group of $A$ is
generated by $\Theta$.
Let $\gamma$ be an effective divisor on $A$ fixed by $K(\Theta)$. Then
there is a positive
integer $a$ such that
$\gamma= a\cdot d_n\cdot \Theta$ in the Neron-Severi group of $A$.}

\vskip0.2truecm\noindent {\it Proof.} Let ${\hat A}$ be the polarized dual
variety of $A$. The
primitive dual polarization ${\hat \Theta}$ is of type $(1,{{d_n}\over
d_{n-1}},...,
{{d_n}\over d_{2}}, d_n)$ and it generates the Neron-Severi group of ${\hat
A}$. We have the
map $\lambda_\Theta: A\to {\hat A}$. Then an easy computation using
self-intersection
numbers shows that:\par

$$\lambda_\Theta^ *({\hat \Theta})=d_n\Theta.$$

\noindent On the other hand we have
$\gamma=\lambda_\Theta^ *({\hat \gamma})$, where ${\hat \gamma}$ is an
effective divisor on
${\hat A}$. Therefore we have ${\hat \gamma}=a{\hat \Theta}$ for some
positive integer $a$. By
pulling this back to $A$ via $\lambda_\Theta$, we get the assertion.
$\diamondsuit$
\vskip0.4truecm

\vskip0.4truecm\noindent {\bf
Lemma} 1.9 {\sl Assume that $\epsilon \not\in K(\Theta)$ and that the
Neron-Severi group of
$A$ is generated by
$\Theta$. Let
$D$ be any irreducible element of the pencil $|-2K_X|$ and let $H_D$ be the
element of the
Neron-Severi group of $D$ given by the restriction of $H$ to $D$.
Then $H_D=\delta_i^ *(\Theta_i)$, $i=1,2,3$ and $H_D$ is at most divisible
by $2$ in the Neron-Severi group of $D$.}

\vskip0.2truecm\noindent {\it Proof.}
The assumption that $NS(A) \simeq {\bf Z}$ implies that also
$NS(D) \simeq {\bf Z}$.
The assertion $H_D=\delta_i^ *(\Theta_i)$ is then purely
numerical and follows  from
the fact that
$D=2A$ in $NS(X)$. Since the maps $\delta_i$ are $2:1$ covers
the maps ${\delta_i}^*: NS(A_i) \simeq {\bf Z} \to NS(D) \simeq {\bf Z}$
have a cokernel which is torsion of order at most $2$.
$\diamondsuit$
\vskip0.4truecm

In the case of abelian surfaces we can extend lemma (1.8) above in the
following way:\par

\vskip0.4truecm\noindent {\bf
Lemma} 1.10 {\sl Let $A$ be an abelian surface with a polarisation
$\Theta$ of type
$(1,2n)$, resp. $(2,n)$ and assume that $NS(A) \simeq {\bf Z}$. If $C$ is a
curve invariant under a group $G \simeq {\bf Z}_n \times {\bf Z}_n$ which
acts on $A$ by translation, then $C=a\cdot \Theta$ where $a$ is a multiple
of $n$, resp. $n/2$.}
\vskip0.2truecm\noindent {\it Proof.} By the assumption $NS(A) \simeq {\bf Z}$
the curve $C$ is a multiple of $\Theta$, resp. $\Theta/2$.
Since $C$ is invariant under the group
$G$ the associated line bundle descends to $A/G$. But this implies that
$G$ is a totally isotropic subgroup with respect to the Weil pairing
(cf [LB, corollary 6.3.5])
of $K(a\Theta)$. This is only possible
if $a$ is divisible by $n$, resp. $n/2$
(cf the description of the Weil pairing given in [LB, example 7.7.4]).
$\diamondsuit$

\vskip 1.0 truecm

\noindent {\bf 2. Some scrolls of secant lines to abelian varieties.}\par

\vskip 0.2truecm

Let us consider a linearly normal abelian variety $A\subset {\bf P}^{r-1}$
of dimension $n$,
embedded via a very ample line bundle ${\cal L}$ belonging to a
polarization $\Theta$ of type
$(d_1,d_2,...,d_n)$. Then $r=d_1\cdot ...\cdot d_n$ and the degree of $A$
equals
$n!\cdot d_1\cdot...\cdot d_n$. Let
$\epsilon$ be a non trivial $2$-torsion point on $A$. We are interested in the
$n+1$-dimensional scroll:\par

$$Y=Y(A,\epsilon, {\cal L})=\cup_{x\in A} L(x,x+\epsilon)$$

\noindent where $L(a,b)$ is the line joining two distinct points $a,b$ in
projective
space. We notice that, unless $n=1$ and $r=d=3$, $Y$ is a proper subvariety
in ${\bf
P}^{r-1}$. \par

As we saw in \S 1, from which we keep the notation, we can associate to
this situation a
${\bf P}^1$-bundle $X$ on ${\overline A}=A/\epsilon$. The relation between
$X$ and $Y$ is
described in the following lemma:\par

\vskip0.4truecm\noindent {\bf
Lemma} 2.1 {\sl One has the following commutative diagram:\par

$$ \xymatrix {
  A \ar[r] \ar[d]_{\phi_{\cal L}} &
    X \ar[d]^{\phi:=\phi_{{\cal O}_X(1)}} \\
  {\bf P}^{r-1}={\bf P}(H^0(A,{\cal L})^*) \ar @{}[r]|{\;\;\simeq} &
    {\bf P}(H^0(X,{\cal O}_X(1))^*)
} $$

%$$ \matrix { & A  & \to   & X \cr
%\phi_{\cal L} & \downarrow  &  &   \downarrow   & \phi:=\phi_{{\cal
%O}_X(1)} \cr
%& {\bf P}^{r-1}={\bf
%P}(H^0(A,{\cal L})^*)   &   \simeq   &   {\bf P}(H^0(X,{\cal O}_X(1))^*)\cr}$$

\noindent Moreover the map $\phi$ is a morphism and its image is $Y$.}

\vskip0.2truecm\noindent {\it Proof.} Since
$p_*{\cal O}_X(1)\simeq {\cal E}\simeq \pi_*{\cal L}$
we have a canonical isomorphism of vector spaces $H^0(A,{\cal L})\simeq
H^0(X,{\cal O}_X(1))$.
Moreover, by lemma (1.2, ii) we have ${\cal O}_X(1)_{|A}\simeq {\cal L}$.
This shows the
existence of the commutative diagram above, which in turn, implies
$\phi(X)=Y$.\par

Let $F$ be any ruling of $X$ and let
$x,x+\epsilon$ be the two points where $F$ intersects $A$.
Since ${\cal L}$ is very ample on $A$, the linear system $|{\cal O}_X(1)|$
separates these
two points and hence $|{\cal O}_X(1)|$ has no base points on $X$, i.e.
$\phi$ is a morphism.
$\diamondsuit$
\vskip0.4truecm

We can now prove the following propositions:\par

\vskip0.4truecm\noindent {\bf
Proposition} 2.2 {\sl The map $\phi$ is finite, i.e. ${\cal O}_X(1)$ is ample.}

\vskip0.2truecm\noindent {\it Proof.} Assume that there is an irreducible
curve $C$ which is contracted under $\phi$. The curve $C$ can only meet a
ruling $F$ once, i.e. the projection of $C$ to $\overline A$ is birational onto
its image. After possibly replacing $C$ by its normalization we obtain a smooth
curve ${\tilde C}$ and a morphism $f:{\tilde C}\rightarrow {\overline A}$ which
is birational onto its image such that
$$
f^{\ast}{\cal E}={\cal O}_{{\tilde C}}\oplus \xi
$$
for some suitable line bundle $\xi$ on ${\tilde{C}}$. The ${\bf P}^1$-bundle
$V={\bf P}(f^{\ast}{\cal E})$ over ${\tilde C}$ is mapped to a cone $W$ in
$Y$. Let $C_0$ be the
section of $V$ which is mapped to the vertex of $W$ and let $C_1$ be the
2-section which is the pullback of the 2-section $A$ of $X$. We denote by $f$
the class of a fibre in $V$. Then $C_1=2 C_0+ af$ in the Neron-Severi group of
$V$ for some integer $a$. Since $C_0$ is contracted it follows that
$e=C_0^2 < 0$. On the other hand it follows from lemma (1.2)(ii) that
$C_1^2=0$ and hence $a=-e$. But then $C_0.C_1=e<0$, i.e. $C_0$ and $C_1$ have
a common component. This contradicts the fact that
the abelian variety $A$ is embedded by the map
$\phi$. $\diamondsuit$
\vskip0.4truecm

\vskip0.4truecm\noindent {\bf
Proposition} 2.3 {\sl One has:\par

$$deg(\phi) \cdot deg(Y)={1\over 2}(n+1)! \cdot r.$$

\noindent If $n=1$ and $r\geq 4$, then $\phi$ is birational onto its
image and it
is an embedding as soon as $r\geq 5$.}

\vskip0.2truecm\noindent {\it Proof.} The result is well known for $n=1$
(see e.g. [CH, proposition 1.1 and proposition 1.2]). So
we assume
$n\geq 2$. We have:

$$deg(\phi) \cdot deg(Y)={\cal O}_X(1)^{n+1}={1\over 2}(f^*{\cal
O}_X(1))^{n+1}={1\over
2}({\cal O}_{\pi^*({\cal E})}(1))^{n+1}.$$

\noindent We can then deduce from formula (1) of \S 1 that:\par

$${\cal O}_X(1)^{n+1}={1\over 2}(n+1)!{\Theta}^n={1\over 2}(n+1)! \cdot r$$
and hence the assertion. $\diamondsuit$
\vskip0.4truecm

In what follows we will need a formula for the double point
cycle ${\bf D}(X,l)$ (see [F, p. 166])
of a map $f: X\to {\bf P}^ l$, where $l\leq r-1$ and $f$ is the
composition of $\phi$
with a general projection ${\bf P}^{r-1}
\to {\bf P}^ l$. We assume $l \geq n+1$. In this situation the map
$f$ is finite. \par

\vskip0.4truecm\noindent {\bf Theorem} 2.4 {\sl One has:
\par

$$ {\bf D}(X,l)= ({1\over 2}(n+1)!\cdot r -{{l+1}\choose {n+2}})
H^{l-n-1}+ {{l+1}\choose {n+3}}H^ {l-n-2}\cdot A$$

\noindent in the homology ring of $X$.}

\vskip0.2truecm\noindent {\it Proof.} By applying theorem (9.3) from [F, p.
166],
one has:\par

$${\bf D}(X,l)= f^ *f_*[X]- c_{l-n-1}(N_f)\cap [X]$$

\noindent where $N_f$ is the normal sheaf to the map $f$ which is defined
by the exact
sequence:\par

$$0\to T_X\to f^ *T_{{\bf P}^ l}\to N_f\to 0.$$

\noindent By proposition (2.3) and finiteness of $f$, we have:

$$f^ *f_*[X]=({1\over 2}(n+1)!\cdot r) H^ {l-n-1}.$$

\noindent Moreover from the Euler sequence, we see that $c(f^ *T_{{\bf P}^
l})=(1+H)^ {l+1}$.
Also by the exact sequence:\par

$$0\to p^*T_{\overline A}\to T_X\to T_{X|{\overline A}}\to 0$$

\noindent we deduce that $c(T_X)=1+A$. Since $A^ 2=0$ in homology, the
assertion follows by an
easy computation. $\diamondsuit$

\vskip0.4truecm \noindent {\bf
Corollary} 2.5 {\sl Let $n\geq 2$ and $r=2n+3$. Then:\par

$${\bf D}(X,r-1) = n!\cdot (2n+3) [{1\over 2} (n+1)
({1\over 2} (n+1)! \cdot (2n+3) - {{2n+3}\choose {n+2}})
+  {{2n+3}\choose {n+3}}]$$

\noindent which is equal to $0$ if and only if $n=2$.}

\vskip0.2truecm\noindent {\it Proof.} The formula for ${\bf D}(X,r-1)$
follows right away
from the above theorem .
For $n=2$ we have ${\bf D}(X,6)=0$ and for $n=3$ one computes that
${\bf D}(X,8)= 6 \cdot 9 \cdot 48$. In order to prove the
second assertion, it is sufficient to show that:\par

$$F(n)={1\over 2}(n+1)! \cdot (2n+3) - {{2n+3}\choose {n+2}}\geq 0$$

\noindent for all $n\geq 4$. Since

$$F(n) > {1\over 2}(n+1)! \cdot (2n+3) - 2^{2n+2}$$
\noindent one can easily verify this by induction on $n$. $\diamondsuit$

\vskip0.4truecm

\vskip0.4truecm \noindent {\bf
Remark} 2.6 The above corollary suggests that the secant scroll to an
abelian surface of type
$(1,7)$ in ${\bf P}^ 6$ should be smooth. This we are going to prove in \S
4. Moreover this
is the only case, apart from the elliptic scroll of degree $5$ in ${\bf P}^
4$, in which an
$(n+1)$-dimensional scroll obtained as above from an $n$-dimensional
abelian variety, which is
linearly normal in ${\bf P}^ {2n+2}$, can be smooth. This is what makes the
consideration of
the surface case, to which the main part of this paper is devoted, particularly
interesting.\par

\vskip0.4truecm We want to finish this section with a remark which is
specific to the case
$\epsilon\in K(\Theta)$.

\vskip0.4truecm \noindent {\bf
Remark} 2.7 If $\epsilon\in K(\Theta)$, then $\epsilon$ acts as an
involution on:

$$H^0(A,{\cal L})\simeq H^0(X,{\cal O}_X(1))
\simeq H^0({\overline A},{\cal E})\simeq H^0({\overline A},
{\overline {\cal L}})\oplus H^0({\overline A},{\overline {\cal L}}\otimes
{\cal M}_1),$$

\noindent the invariant and anti-invariant eigenspaces being
$H^0({\overline A},
{\overline {\cal L}})$ and  $H^0({\overline A},{\overline {\cal L}}\otimes
{\cal M}_1)$
respectively.  Recall that $\overline {\cal L}$ represents a polarization
of type
$(\overline {d_1},...,\overline {d_n})$ on
$\overline A$, such that $2\cdot \overline {d_1}\cdot ...\cdot \overline
{d_n}={d_1}\cdot
...\cdot {d_n}$. Set $\overline r=\overline {d_1}\cdot ...\cdot \overline
{d_n}$. Then
$h^0({\overline
A},{\overline {\cal L}})=h^0({\overline A},{\overline {\cal L}}\otimes {\cal
M}_1)={\overline r}$ and
$2{\overline r}=r$.\par

 Accordingly $\epsilon$ acts as an involution on ${\bf P} ^{r-1}={\bf
P}(H^0(A,{\cal
L})^*)$. The invariant and anti-invariant subspaces both have dimension
${\overline r}$ and
are ${\bf P}^+={\bf P}(H^0({\overline A},{\overline {\cal L}})^*)$ and
${\bf P}^-={\bf P}(H^0({\overline A},{\overline {\cal L}}\otimes {\cal
M}_1)^*)$. One has the morphisms $\phi^+=\phi_{\overline {\cal L}}:
{\overline A}\to {\bf
P}^+$ and $\phi^-=\phi_{\overline {\cal L}\otimes {\cal
M}_1}: {\overline A}\to {\bf
P}^-$. The images $Y^+$ and $Y^-$ of these maps are nothing but the images
via $\phi$ of
the two sections ${\overline A}^+$ and ${\overline A}^-$ respectively of
$X$ (see remark
(1.4, i)). \par

Finally we have a different description of $Y$ which will be useful to take
into account.
Take any point $x\in {\overline A}$ and consider the corresponding points
$x_{\pm}\in {\overline A}^{\pm}$. Set $y^{\pm}=\phi(x^{\pm})$. Then:

$$Y=\bigcup_{x\in {\overline A}}
L(y^+,y^-).$$

\vskip 2.0 truecm

\noindent {\bf 3. Secant scrolls related to abelian surfaces.}\par

\vskip 0.2truecm

From now on we will consider the case where the abelian variety $A$
is a surface, i.e. $n=2$,
and the polarization $\Theta$ is very ample of type $(d_1,d_2)$.
Hence $r=d_1\cdot d_2$ and the
surface
$A$ is embedded into ${\bf P}^ {r-1}$, via a line bundle ${\cal L}$
representing
$\Theta$, as a
surface of degree $2\cdot d_1\cdot d_2$. We will also denote the image by
$A$. \par

We are interested in characterizing the cases in which the map $\phi: X\to
Y$ introduced in
general in
\S 2 is an embedding, if $A$ is a general polarized abelian surface of type
$(d_1,d_2)$.
Notice that, by Lefschetz's hyperplane section theorem, there is
no chance that
$\phi$ is ever an embedding if $r<7$ .
So we will assume from now on $r\geq 7$. \par

In order to study the map $\phi: X\to Y$, we need some information about
the embedding of $A$ in
${\bf P}^ {r-1}$. We will use a well known result of Reider, in the
following form due to
Beltrametti and Sommese [BS, theorem 3.2.1]:

\vskip0.4truecm\noindent {\bf
Theorem} 3.1 {\sl  Let $L$ be a
numerically effective (nef) divisor on a surface $S$. Assume that $L^2\geq
4k+1$. Given any
0-dimensional scheme $Z$ of length $k$ on $S$, then either the natural
restriction map

$$H^0(S,{\cal O}_S(K+L))\to H^0(S,{\cal O}_Z(K+L))$$

\noindent is surjective, or there exist an effective divisor $C$ on $S$ and
a non-empty
subscheme $Z'$ of $Z$ of length $k'\leq k$, such that:\par

\noindent (i) the map

$$H^0(S,{\cal O}_S(K+L))\to H^0(S,{\cal O}_{Z'}(K+L))$$

\noindent is not surjective;\par

\noindent (ii) $Z'$ is contained in $C$ and there is an integer $m$ such
that $m(L-2C)$ is effective;\par

\noindent (iii) one has

$$L\cdot C-k'\leq C^2<{{L\cdot C}\over 2}<k'.$$}\par

\vskip0.4truecm

As a consequence we have the following:

\vskip0.4truecm\noindent {\bf
Proposition} 3.2 {\sl Let $A$ be a polarized abelian surface of type
$(d_1,d_2)$ and
let
${\cal L}$ be a line bundle on $A$ representing the given polarization
$\Theta$.
Set $r=d_1\cdot d_2$. One has:\par

\noindent (i) if  $r\geq 5$, then ${\cal L}$ is very ample, unless there is
a curve
$C$ on $A$ such that $C^2=0$ and $\Theta\cdot C\leq 2$;\par

\noindent (ii) assume that  ${\cal L}$ is very ample and that $r\geq 7$,
then $\phi_ {\cal
L}(A)$ has no
$3$-secant lines, unless there is a curve
$C$ on $A$ such that either $C^2=0$ and $\Theta\cdot C=3$ (a plane cubic)
or $C^2=2$ and
$\Theta\cdot C= 5$ (a genus $2$ quintic);\par

\noindent (iii) assume that  ${\cal L}$ is very ample, that $\phi_ {\cal
L}(A)$ has no
$3$-secant lines and that $r\geq 9$, then $\phi_ {\cal
L}(A)$ has no $4$-secant planes, unless there is a curve
$C$ on $A$ such that either $C^2=0$ and $\Theta\cdot C=4$ (a genus $1$
quartic) or $C^2=2$
and $\Theta\cdot C=6$ (a genus $2$ sextic).\par

\noindent In particular if $NS(A)\simeq {\bf Z}$
then:\par

\noindent (i') if  $r\geq 5$, then ${\cal L}$ is very ample;\par

\noindent (ii') if $r\geq 7$, then $\phi_ {\cal
L}(A)$ has no
$3$-secant lines;\par

\noindent (iii') if $r\geq 9$, then $\phi_ {\cal
L}(A)$ has no $4$-secant planes unless $\Theta$ is the triple of a principal
polarization.\par }\par

\vskip0.2truecm\noindent {\it Proof.} The first part is an immediate
application of theorem
(3.1), by taking into account that for any effective divisor $C$ on $A$,
the integer $C^2$
is even and non-negative. \par

If $NS(A)\simeq {\bf Z}$, then there is no curve on $A$ with $C^2=0$. In
addition, if $C$ is
an effective divisor such that $C^2=2$, then $C$ is irreducible and its
class $\theta$ in
$NS(A)$ is a principal polarization which is indivisible, hence it
generates $NS(A)$. Thus
there is a positive integer $a$ such that $\Theta=a\theta$. If $\Theta\cdot
\theta=2a\leq
6$, then $a\leq 3$, and $\Theta$ is a polarization of type $(a,a)$. Since
we are
assuming $r=a^2\geq 5$, we have $a=3$, a case which indeed gives rise to
$4$-secant
planes (see \S 6). $\diamondsuit$

\vskip0.4truecm

As a further consequence we can prove the:

\vskip0.4truecm\noindent {\bf
Theorem} 3.3 {\sl Let $A\subset {\bf P}^{r-1}$ be a linearly normal, smooth
abelian surface such that ${\cal L}={\cal O}_A(1)$ determines a polarization
$\Theta$ of type
$(d_1,d_2)$ with $r=d_1\cdot d_2$. Let $\epsilon\in A$ be a non trivial
point of
order two. Suppose that $NS(A)\simeq {\bf Z}$. Consider the map $\phi: X\to
Y\subset
{\bf P}^{r-1}$. Then:\par

\noindent (i) distinct rulings of $X$ are sent by $\phi$ to
distinct lines in ${\bf P }^ {r-1}$, as
soon as $r\geq 7$;\par

\noindent (ii) the differential of the map $\phi: X\to Y$ is injective along
$A$ and $Y$ is smooth along $A$, as soon as
$r\geq 7$;\par

\noindent (iii) the map $\phi: X\to Y$ is an
isomorphism as soon as $r\geq 9$, unless $\Theta$ is the triple of a principal
polarization, i.e. $d_1=d_2=3$ and $r=9$.}

\vskip0.2truecm\noindent {\it Proof.} (i) If two distinct rulings of $X$ were
mapped to the same
line $L$ in ${\bf P }^ {r-1}$, then $L$ would be a $4$-secant line to $A$,
contradicting
corollary (3.2, ii').\par

(ii) Let $x$ be a point of $A$. Then the tangent space to $X$ at $x$ is
spanned by $T_{A,x}$
and by the tangent space to the ruling $F$ through $x$. Since $L=\phi(F)$
cannot be tangent to
$A$ by (3.2, ii'), it follows that $d\phi$ is injective at $x$.
The smoothness of $Y$ along $A$ is
then a consequence of (3.2, ii'). \par

(iii) If $\phi$ were not injective, we would
have two distinct secant lines to $A$ meeting at a point, and therefore a
$4$-secant plane to $A$, contradicting (3.2, iii').\par

Suppose $d\phi$ is not injective at a point $z\in X$, which we may
assume not to be on $A$ by (ii).  The
gaussian map $\gamma: X\to {\bf G}(3,r-1)$ is a rational map whose restriction
to each ruling is defined by (ii). By following the argument in [R, p. 215], we
see that such a restriction is given by quadratic forms. If $F$ is the
ruling of
$X$ through $z$, the forms defining $\gamma_{|F}$ all vanish at $z$. Hence by
[R, lemma 25], the union of the tangent spaces to $Y$ along the image $L$
of $F$
is a ${\bf P}^4$. In particular the two tangent planes to $A$ at the
points $x$ and
$x+\epsilon$ where $L$ intersects $A$ meet at a point. This either yields
the existence of
a tangent line to $A$ which meets $A$ once more, which is impossible, or of two
tangent lines $r$, $r'$ to $A$ at  $x$ and $x'=x+\epsilon$ which meet,
hence the
existence of a
$4$-secant plane to
$A$, a contradiction to (3,2, ii', iii'), unless $d_1=d_2=3$, $r=9$ and
$\Theta=3\theta$, with $\theta$ a principal polarization on $A$. $\diamondsuit$
\vskip0.4truecm

We also have the following consequence:\par

\vskip0.4truecm\noindent {\bf
Proposition} 3.4 {\sl Let $A$ be
a polarized abelian surface of type $(d_1,d_2)$ with $r=d_1\cdot
d_2\geq 7$ such that $NS(A)\simeq {\bf Z}$.
Then the map
$\phi: X\to Y$ is birational.}

\vskip0.2truecm\noindent {\it Proof.} Let $\delta$ be the degree of $\phi$.
Let $F$ be a general
ruling of $X$ and let $L$ be its image under $\phi$.
Let $x$, $x+\epsilon$ be the two points of
$A$ on $L$. Let
$z\in L$ be a general point. Then by (3.2, ii')
there are $\delta-1$ images of rulings of $X$ through
$z$ different from $L$. This situation produces a $(\delta-1)$-cover
of $L$ and again by (3.2, ii') this
cover is totally ramified at $x$ and $x+\epsilon$,
i.e. ${\phi}^{-1}(\phi(x))=x$ and
${\phi}^{-1}(\phi(x+\epsilon))=x+\epsilon$. Since $z$ is general this implies
that $\phi$ itself should
be ramified along the points of $A$, contradicting (3.3, ii).
$\diamondsuit$
\vskip0.4truecm

In what follows we will need some information about the hyperplane
sections of $Y$. Let $A\subset {\bf P}^{r-1}$ be a linearly normal, smooth
abelian surface such that ${\cal L}={\cal O}_A(1)$ determines a polarization
$\Theta$ of type
$(d_1,d_2)$ with $r=d_1\cdot d_2$. Let $\epsilon\in A$ be a non trivial
point of
order two. Let
$H$ be a divisor in $|{\cal O}_X(1)|$. We abuse notation and we
denote by $p: H\to
{\overline A}$ the restriction to $H$ of the projection $p: X\to {\overline
A}$. \par

\vskip0.4truecm\noindent {\bf
Lemma} 3.5 {\sl In the above setting, if $NS(A)\simeq {\bf Z}$
 generated by $\Theta$
 and if
$H$ is a general divisor in $|{\cal O}_X(1)|$, then:\par

\noindent (i) the map $p: H\to {\overline
A}$ is the blow-up of ${\overline A}$ at $r$ distinct points $p_1,...,p_r$
and the
exceptional divisors are $r$ rulings $F_1,...,F_r$ of $X$ contained in $H$;\par

\noindent (ii) if $\epsilon \in K(\Theta)$,
then ${\cal O}_X(1)_{|H}= {\cal O}_X(-K_X)\otimes {\cal
O}_H(F_1+...+F_r)={\cal O}_X(A) \otimes {\cal O}_H(F_1+...+F_r)$, whereas
if $\epsilon \notin K(\Theta)$, then ${\cal O}_X(1)_{|H}=
{\cal O}_X(-K_X)\otimes {\cal O}_H(F_1+...+F_r)={\cal O}_X(A)\otimes {\cal
M}_1
\otimes {\cal
O}_H(F_1+...+F_r)$. }

\vskip0.2truecm\noindent {\it Proof.} Let $H$ be the zero locus of the
section $s\in H^ 0(X,
{\cal O}_X(1))$. One has $H^ 0(X,
{\cal O}_X(1))\simeq H^ 0(A,{\cal L})$ and  $H^ 0(X,
{\cal O}_X(1))\simeq H^0(\overline A,{\cal E})$, thus we may
interpret $s$ as a non-zero section of ${\cal L}$ on $A$ and of ${\cal E}$ on
$\overline A$. \par

We let $C_H$ be the zero-divisor of $s\in H^ 0(A,{\cal L})$ on $A$, i.e.
the intersection of $H$ with $A$. According to our assumption on $A$,
 the curve $C_H$ is irreducible and reduced. Then we consider the curve
$C_{H,\epsilon}=C_{H+\epsilon}$,
which is the zero locus
of the section $s_{H,\epsilon}=t_\epsilon^ *(s)\in H^ 0(A,t_\epsilon^ *
{\cal L})$. The curve $C_{H,\epsilon}$ is also irreducible and reduced.
Furthermore, it is distinct from $C_H$. This is clear if $\epsilon\notin
K(\Theta)$,
whereas, if $\epsilon\in K(\Theta)$ it follows from the considerations in
remark (2.7).
Hence we can consider the zero-cycle $Z$ of lenght $2r$ given by the
intersection of
$C_H$ with
$C_{H,\epsilon}$. \par

Let $\overline Z$ be the zero locus of $s\in H^0(\overline A,{\cal E})$ on
$\overline A$. It is clear that $\pi^*(\overline Z)=Z$. Hence the
general section of ${\cal E}$ vanishes in codimension two. Notice that this
fits with the fact that $c_2({\cal E})={{\Theta^2}\over2}=r$.

Now we claim that, by the
generality assumption we are making on $H$,
the cycle $\overline Z$ consists of $r$ distinct points. In view of the above
considerations, this amounts to proving that a general section $s\in
H^0(\overline A,{\cal
E})$ vanishes at $r$ distinct points of $\overline A$. This is a
consequence of the
following claim, which is a Bertini type theorem for vector bundles:

\vskip0.2truecm\noindent {\it Claim}: {\it Let $V$ be a smooth irreducible
variety and let
${\cal E}$ be a vector bundle of rank two on $V$. Assume that:\par

\noindent (i) $NS(V)\simeq {\bf Z}$
is generated by $c_1({\cal E})$,\par

\noindent (ii) ${\cal E}$ is generated by global sections away from a
subvariety
of $V$ of codimension $\geq 3$.\par

\noindent Then the zero locus of the general section of ${\cal E}$ is
reduced.}\vskip0.2truecm

As for the proof, let $\sigma, \tau$ be general sections of ${\cal E}$. Then
$\sigma\wedge\tau$ is not identically zero by (ii).
Let $D$ be the zero locus of
$\sigma\wedge\tau$. By (i), the subvariety $D$ is an
irreducible, reduced divisor on
$V$. In particular this implies that the general section of ${\cal E}$
vanishes in
codimension $2$ and, again by (ii), the zero loci of two general sections,
like
$\sigma, \tau$ of ${\cal
E}$ have no common component. \par

Set $s(\lambda,\mu)=\lambda\sigma+\mu \tau$, with $[\lambda,
\mu]\in {\bf P}^1$. The zero locus $W(\lambda,\mu)$ of $s(\lambda,\mu)$
describes, as
$[\lambda,
\mu]$ varies in ${\bf P}^1$, a linear system of divisors on $D$, which, by
the above
considerations, has no fixed divisor. Hence Bertini's theorem ensures that
for $[\lambda,
\mu]$ general in
${\bf P}^1$, the scheme
$W(\lambda,\mu)$ is reduced. \par\vskip0.2truecm

Let us now return to $\overline Z$.
By the above claim (Note that condition (ii) is fulfilled by
lemma (2.1)) we can write ${\overline
Z}=p_1+...+p_r$, with
$p_1,...,p_r$ distinct points of $\overline A$. Then
$H$ contains the rulings
$F_1,...,F_r$ over the points
$p_1,...,p_r$. The map $p: H-\sum_{i=1}^ rF_i\to {\overline A}-Z$ is
clearly an isomorphism
and since the $F_i$'s are contracted by $p$ to the smooth points $p_i$, it
follows that $p:
H\to {\overline A}$ is the blow-up in ${\overline Z}$. This proves part
(i).\par

As for part (ii), we apply the adjunction formula and obtain:\par

$${\cal O}_H(F_1+...+F_r)={\cal O}_H(K_H)=
{\cal O}_H\otimes {\cal O}_X(1)\otimes {\cal O}_X(K_X)$$

\noindent which, by proposition (1.6), concludes the proof. $\diamondsuit$

\vskip 1.0 truecm

\noindent {\bf 4. The case $r=7$.}\par

\vskip 0.2truecm

In this section we will prove the:\par

\vskip0.4truecm\noindent {\bf
Theorem} 4.1 {\sl Let $A$ be
a linearly normal abelian surface embedded in ${\bf P}^6$ via a line bundle
${\cal L}$ which
determines a polarization $\Theta$ of type $(1,7)$. Assume that
$End(A)\simeq {\bf Z}$. Then the map $\phi: X\to Y$ is an
embedding.}

\vskip0.4truecm Fist of all we notice
that
$End(A)\simeq {\bf Z}$ implies
$NS(A)\simeq {\bf Z}$ (see [LB, p. 122]). Since $\Theta$ is indivisible, it
generates $NS(A)$. Next we remark that no non-trivial $\epsilon$ is an element
of $K(\Theta)$.\par

 The proof will require several steps. Let us start by denoting by
$\Sigma$ the singular locus of $Y$
and let $S=\phi^ {-1}(\Sigma)$. Recall that $|2A|$ is a base point free pencil
on $X$ containing
three double fibres $A=A_1$, $A_2$, and $A_3$, which are the only
reducible, singular members of the
pencil (see \S 1, from which we keep the notation). Notice also that
$End(A_i)\simeq {\bf Z}$, for $i=1,2,3$, since we have made the assumption that
$End(A)\simeq {\bf Z}$.
In particular
${\cal L}_i$ is very
ample on $A_i$, $i=1,2,3$, and $A_1$, $A_2$, $A_3$
play a symmetric role with respect to $X$ and $Y$.\par

\vskip0.4truecm\noindent {\bf
Lemma} 4.2 {\sl One has:

$$ S\cap A_i=\emptyset,\quad i=1,2,3 \eqno (1)$$

\noindent and for every irreducible
component $Z$ of $S$, there is a unique irreducible, smooth surface
$D\in |2A|$ on $X$ containing $Z$.}

\vskip0.2truecm\noindent {\it Proof.} As we noticed, $A_1$, $A_2$, $A_3$
play a symmetric role with respect to $X$ and $Y$. Therefore (1) follows by
(3.3, ii).
Moreover, $Z\cdot A=0$ in the homology ring and $|2A|$ is a base point free
pencil,
whose elements sweep out $X$. Therefore there is an element $D\in |2A|$
containing
$Z$. By proposition (1.5) the surface $D$ is smooth and irreducible.
$\diamondsuit$
\vskip0.4truecm

We want to prove that $\Sigma=\emptyset$. First we prove that:\par

\vskip0.4truecm\noindent {\bf
Lemma} 4.3 {\sl One has $dim(\Sigma)\leq 1$.}

\vskip0.2truecm\noindent {\it Proof.} We argue by contradiction and
therefore, according to lemma
(4.2), we may assume that one of the following happens:\par

\noindent (i) there is a surface $D\in |2A|$ on which $d\phi$ is not
injective;\par

\noindent (ii) there are surfaces $D_1, D_2\in |2A|$ which have the same
image via $\phi$;\par

\noindent (iii) there is a surface $D\in |2A|$ on which $\phi$ is not
injective.\par

In case (i) the differential $d\phi$ would not be injective at $4$ distinct
points of a
general ruling $F$ of $X$. By the same argument we made in the proof of
part (iii) of theorem
(3.3), $d\phi$ would not be injective along the whole of $F$, a
contradiction.\par

In case (ii) we may assume that both $D_1$ and $D_2$ map birationally, via
$\phi$, to some
irreducible component of $\Sigma$, otherwise we are in case (iii). Hence we
have a birational
map $a: D_1\to D_2$, which is an isomorphism, since $D_1$ and $D_2$ are both
abelian surfaces.
Recall that $D_1$ and $D_2$ are isomorphic as double covers of $A$ and that
this isomorphism is compatible with the projection onto $\overline A$.
Since $End(A)\simeq {\bf Z}$ the same is true for $D_1$ and $D_2$.
Hence  the map $a$ is of type $a: z\in D_1\to \pm z+ k\in D_2$, where
$k$ is a fixed
element in $D_2$. But then this would imply that the $4$ points of intersection
of $D_1$ with any
ruling $F$ of $X$ are mapped, via $a$, to the $4$ points of intersection of
$D_2$ with another
ruling $F'$. Hence $\phi$ would map $F$ and $F'$ to the same line in ${\bf
P}^ 6$,
contradicting theorem (3.3, i).
\par

In case (iii), let $\mu$ be the degree of $\phi_{|D}$. Since $H^ 2\cdot
D=28$, we can
only have the possibilities $\mu=2,4,7,14,28$. If  $\mu\geq 7$, then
$\Sigma'=\phi(D)$ would be degenerate, implying that the linear system
$|H-D|=|H-2A|$ is
effective, a contradiction, since $F\cdot (H-2A)=-3$. \par

If $\mu=2$, then we have an involution $a$ on $D$ and we can argue as we
did in case (ii) to
arrive at a contradiction.\par

Let us consider the last case $\mu=4$, in which case the degree of $\Sigma'$ is
$7$, and we may
assume that $\Sigma'$ is non-degenerate in ${\bf P}^ 6$. Let $z$ be a
general point of $D$ and
let $F$ be the ruling of $X$ through $z$.
Since $\phi$ maps $F$
isomorphically to a line $L$
in  ${\bf P}^ 6$, every such ruling is a proper $4$-secant to the
surface $\Sigma'$, i.e. it intersects the surface in $4$ different points.
We see that the fibres containing $z$ of the maps $D\to
{\overline A}$ and $D\to
\Sigma'$ intersect only at $z$. This implies that through the general point
$p=\phi(z)\in
\Sigma'$ there are at least four $4$-secant lines to $\Sigma '$. Hence, by
projecting
$\Sigma'$ from $p$ down to ${\bf P}^ 5$, we have a surface $\Sigma''$ of
degree $6$ with at
least $4$ distinct triple points $p_1,p_2,p_3,p_4$. \par

First we claim that $p_1,p_2,p_3,p_4$ are not collinear for
a general projection. Suppose in fact
they lie on a line
$L$. Then a simple application of Bezout's theorem shows that the general
hyperplane section
$\Gamma$ of $\Sigma''$ through $L$ consists of $2L$ plus a residual curve
$\Gamma'$
of degree $4$ containing $p_1,p_2,p_3,p_4$. Hence, again by Bezout's
theorem, $\Gamma'$
must be
reducible into degenerate components. This implies that the projection of
$\Sigma''$ from $L$
to ${\bf P}^ 3$ is a non-degenerate curve $\Lambda$. Let $\alpha\geq 2$ be
the multiplicity
of
$\Sigma''$ at a general point of $L$, let $\lambda\geq 3$ be the degree of
$\Lambda$ and let
$\beta$ be the degree of the general fibre of the projection $\Sigma''\to
\Lambda$. Then we
have:\par

$$6=\alpha+\lambda\beta\geq 2+\lambda\beta.$$

\noindent This gives us the possibilities: $\alpha=2$,
$\lambda=4$ and  $\beta=1$ or $\alpha=3$, $\lambda=3$ and $\beta=1$.
In either case $\Sigma''$
would be a scroll, and $\Gamma'$ would
consist of $3$ or $4$ rulings.
We shall first treat the case where $\alpha=2$.
Since $\Sigma''$ has a double
line we have two possibilities. Either $\Sigma'$ has a double line or
$\Sigma'$ has a pencil of plane curves (recall that the general projection
of $\Sigma'$ contains a double line by our assumption). In the second case
the surface $\Sigma'$ cannot span ${\bf P}^ 6$. This follows since in
this case every ruling meets each of these plane curves in one point (recall
that $\beta=1$) and this shows that all rulings are incident to two planes,
i.e. contained in a fixed ${\bf P}^ 5$. Finally assume that $\Sigma'$ has
a double line. Recall that there are $4$ proper $4$-secant lines through
a general point
of this line. Hence the multiplicity of the corresponding
singular points is greater than the multiplicity of a general point of
the line $L$
and we see that
when $\Gamma'$ moves, there are
infinitely many rulings
through each of the points $p_1,p_2,p_3,p_4$. Hence $\Sigma'$ would be a
cone with vertex each
one of the points $p_1,p_2,p_3,p_4$, a contradiction.
It remains to treat the case where $\alpha=3$. If $\Sigma'$ does not have a
multiple line then we can argue as above that it cannot span ${\bf P}^ 6$.
Finally assume that $\Sigma'$ has
a multiple line. Again, since there are $4$ proper $4$-secant lines through
a general point
of this line, the multiplicity of the corresponding
singular points is greater than the multiplicity of a general point of
the line $L$ and our above argument goes through
unchanged.
\par

Let us now project $\Sigma''$ to ${\bf P}^ 4$ from one of the points
$p_1,p_2,p_3,p_4$. If the
projection is not a surface, then $\Sigma''$ is a cone with vertex the
point we are projecting
from. This cannot happen for all the points under consideration, therefore
we can assume that
the projection from, say, $p_1$ is a surface, which is an irreducible,
non-degenerate
surface of degree $3$ in
${\bf P}^ 4$, with at least two points of multiplicity at least three, a
contradiction.
$\diamondsuit$\par

\vskip0.4truecm

The next step is to further bound the dimension of the singular locus
$\Sigma$ of $X$. \par

\vskip0.4truecm\noindent {\bf
Lemma} 4.4 {\sl One has $dim(\Sigma)\leq 0$.}

\vskip0.2truecm\noindent {\it Proof.} Again we argue by contradiction and
we assume that the locus
$\Sigma$ has components of dimension $1$. Hence $S=\phi^ {-1}(\Sigma)$ has
components of
dimension $1$. Let $S_1$ be their union. According to lemma (4.2), the
curve $S_1$ is
contained in the union of finitely many irreducible divisors of $|2A|$. Let
$D$ be one of
these and let $\gamma$ be the part of $S_1$ contained in $D$. By lemma
(1.9) the polarization $H_D$ is of type $(1,14)$ and hence there is a
positive integer $a$ such that
$\gamma=aH_D$ in $NS(D)$.\par

Since the group $K(\Theta)$ fixes ${\cal L}$, it also acts on the vector
bundle ${\cal E}$
and therefore it acts on $X$. Furthermore it acts on $H^ 0(X, {\cal
O}_X(1))\simeq H^
0(A,{\cal L})$. This action is trivial on the pencil $|-K_X|=|2A|$, since
the divisors $A_i,
i=1,2,3,$ are fixed. Notice that $S$ is fixed by this action, hence $S_1$
is and therefore
$\gamma$ is. \par

Consider now a general map $f: X\to {\bf P}^ 5$ as in \S 2. A
straightforward parameter count
shows that the singular locus of $f(X)$ has still dimension $1$. Then by
applying theorem
(2.4), we see that ${\bf D}(X,5)=6H\cdot (H+A)$. Hence $6H\cdot
(H+A)-\gamma$ is represented by an effective cycle. By lemma (1.10) we have
$\gamma=7aH_D$ in
$NS(D)$, with $a$ a positive integer. Hence we have $\gamma=14aH\cdot A$.
Therefore
$6H^ 2+(6-14a)H\cdot A$ is represented by an effective cycle. By lemma
(3.5, ii), we have $H^
2=H\cdot A+7F$. Hence $(12-14a)H\cdot A+42F$ is represented by an effective
cycle, whose
intersection with the pull-back $W$ on $X$ of an ample divisor on
${\overline A}$ is non
negative. This forces $a\leq 0$, a contradiction. $\diamondsuit$

\vskip0.4truecm

Finally we can finish the:\par

\vskip0.2truecm\noindent {\it Proof of theorem (4.1).} Since we know now
that $S=\phi^
{-1}(\Sigma)$ is finite, we can use the double point formula from corollary
(2.5), which tells
us that ${\bf D}(X,6)=0$, implying that $S$ is empty. $\diamondsuit$

\vskip 2.0 truecm

\noindent {\bf 5. The case $r=8$.}\par

\vskip 0.2truecm

In the case $r=8$ we have two possibilities, namely the abelian surface $A$
is either embedded
in ${\bf P}^7$ via a line bundle ${\cal L}$ belonging to a polarization
$\Theta$ of
type $(1,8)$ or to a polarization $\Theta$ of type $(2,4)$. In the former
case there are three non
trivial points
$\epsilon\in A$ of order two which are contained in $K(\Theta)\simeq {\bf Z}/8
\times {\bf Z}/8$, in the latter case every such
point is an element of
$K(\Theta)$.\par

As we will see in a moment, the two cases $\epsilon\notin K(\Theta)$ and
$\epsilon\in K(\Theta)$ give rise to a
completely different behaviour
of the map
$\phi: X\to Y$.

The main result of this section is the following: \par

\vskip0.4truecm\noindent {\bf
Theorem} 5.1 {\sl Let $A$ be an abelian surface such
that $End(A)\simeq {\bf Z}$, linearly normally
embedded in ${\bf P}^7$ via a line bundle ${\cal L}$ giving a
polarization $\Theta$ in $NS(A)$. Let $\epsilon$ be a non-trivial
 point of order two on
$A$. Then:\par

\noindent (i) if $\epsilon\in K(\Theta)$, the map $\phi: X\to
Y$ fails to be an embedding along
the sections ${\overline A}^\pm$ of $X$ (see remark (1.4, i)), whereas
it is an embedding
on the open subset which is the complement of these two sections;\par

\noindent (ii) if $\epsilon\notin K(\Theta)$, which implies $\Theta$ to be
 of type $(1,8)$,
the map
$\phi: X\to Y$ is an embedding.}

\vskip0.2truecm\noindent {\it Proof.} (i) We use the notation of remark (2.7).
The spaces ${\bf P}^{\pm}$ are both of dimension $3$ and the
maps $\phi^{\pm} :{\overline A}^{\pm}\to Y^{\pm}$
cannot be embeddings. On the other hand $\phi$ is an embedding on
$X-({\overline A}^+\cup
{\overline A}^-)$. This follows from theorem (3.3, ii): the linear system
$|A|$
is a base point free pencil whose elements, with the exception of
$2{\overline A}^+$ and $2{\overline A}^-$, are smooth and isomorphic to
$A$. The pencil $|A|$
sweeps out $X$ and its smooth elements
play a symmetric role in the description of $X$ and $\phi$.\par

(ii)  The structure of the proof in this case is somewhat similar to the one of
the case $d_1=1, d_2=7$, $r=7$. We still denote by $\Sigma$ the
singular locus of $Y$ and we let $S=\phi^ {-1}(\Sigma)$. Lemma (4.2) still
holds. The proof of
lemma (4.3) can be adapted with minor changes to the present situation,
showing again that the
dimension of $\Sigma$ is at most $1$.
We leave these details to the reader.  \par

Let now $S_1$ be the union of the one-dimensional components of $S$. Again
$S_1$ is
contained in the union of finitely many irreducible divisors $D$ of $|2A|$.
By theorem (2.4) we find for the double locus ${\bf
D}(X,5)=9H^2+6H\cdot A=15H\cdot A+72F$ where the last
equality follows from theorem
(3.5). By lemma (1.9) the polarization $H_D$ is of type $(1,16)$ or
$(2,8)$.
By lemma (1.10) there is an integer $a$ such that $S_1=4a H\cdot D=8a
H\cdot A$ in homology. Since ${\bf D}(X,5)-S_1=(15-8a) H\cdot A+72F$ is
effective it follows that $a=1$ and that $S_1$ is contained in a unique surface
$D\in|2A|$.\par
We first remark that $S_1$ is reduced by lemma (1.10). It then follows that the
map $\phi$ restricted to $S_1$ is generically two-to-one onto
the curve $\Gamma=\phi(S_1)$.
(Here we use that if a point $P$ is simple for the double point cycle, then
the differential at $P$ is injective and there is only one other point $Q$
mapped to the same point as $P$. This can be deduced from the construction of
the double point locus (cf.[F, p.166]).) But then the degree of $\Gamma$ is
$64$. Let $C$ be a general tangent hyperplane section of $\phi(D)$. This is an
irreducible reduced curve of degree $32$ with $65$ nodes. Its pullback on $C$
has self-intersection $32$ and one node, hence geometric genus $16$. Thus the
arithmetic genus of $C$ is $81$ which is equal to the Castelnuovo bound for
non-degenerate curves in ${\bf P}^6$ (cf. [EH, p. 87]). Since this bound can
only be achieved by smooth curves we obtain a contradiction.
\par

We may therefore assume that the singular locus $\Sigma$ of $Y$ is
finite. Assume $\Sigma$ is not empty. Then $S=\phi^ {-1}(\Sigma)$
consists of orbits of
$K(\Theta)\simeq {\bf Z}_8^ 2$. From theorem (2.4) we see that ${\bf
D}(X,6)=72$. This
implies that $S$ consists of a unique orbit formed by $64$ distinct points
that are
pairwise coupled by $\phi$, which sends them to $32$ distinct points of $Y$
where two
branches of
$Y$ meet transversally. Let $Z$ be this set of $32$ points in ${\bf P}^ 7$.
\par

Since
every surface
$D\in |2A|$ is
$K(\Theta)$-invariant, $S$ must lie on a unique surface $D\in |2A|$, which
is therefore
mapped by $\phi$ to a surface $\Delta$ with at least $32$ double points at
$Z$.
Since $Z$ is also
$K(\Theta)$-invariant, there is for each $z\in Z$ an element $h_z$ of order
$2$ of $K(\Theta)$
fixing it. The element $h_z$ does not depend on $z$. This follows since the
stabilizers of elements in the same orbit for a group action are conjugated,
resp. equal if the group is abelian.

Now there are three elements of order two in $K(\Theta)\simeq {\bf Z}_8^
2$. Remember how
$K(\Theta)$ acts on ${\bf P}^ 7$ where the coordinates are $[x_0,...,x_7]$.
We have the two
generators $\sigma$ and $\tau$ of $K(\Theta)$ acting as follows
(cf. [LB, p. 169]):\par

$$\sigma: [x_0,...,x_7]\to [x_7,x_0...,x_6,]$$
$$\tau: [x_0,...,x_7]\to [x_0,\xi^ {-1} x_1,..., \xi^ {-i} x_i,..., \xi^
{-7} x_7]$$

\noindent where $\xi=exp({{2\pi\sqrt{-1}}\over 8})$. The elements of order
two are $\sigma^
4$, $\tau^ 4$ and their product. Suppose $h=\tau^ 4$ (the discussion is
similar in the
other cases). Then: \par

$$h: [x_0,...,x_7]\to [x_0, -x_1,..., (-1)^{i}x_i,..., -x_7]$$

\noindent and therefore its eigenspaces are both of dimension $3$. Since
$Z$ is contained in
their union, we deduce that at least $16$ points of $Z$ lie in a ${\bf P}^
3$ which we
denote by $P$. Notice that $P\cap \Delta$ is finite. Otherwise
we could deduce by applying $\sigma$ that both
eigenspaces of $h$ would intersect $\Delta$
in a curve. But then $D$ would contain $2$ curves which do not intersect
which contradicts our assumption that $NS(A)\simeq {\bf Z}$ and hence also
$NS(D)\simeq {\bf Z}$.
Now consider the
intersection of
$\Delta$ with a hyperplane through $P$. This is an irreducible curve $B$ of
degree
$32$, non degenerate in ${\bf P}^ 6$, with at least $16$ singular points on
$P$. Let $q$ be
a point on $B$. Any hyperplane in ${\bf P}^ 6$ containing $P$ and $q$ has
to contain $B$ by
Bezout's theorem, a contradiction. $\diamondsuit$

\vskip 1.0 truecm

\noindent {\bf 6. The case $r=9$.}\par

\vskip 0.2truecm

In this section we prove the following\par

\vskip0.4truecm\noindent {\bf
Theorem} 6.1 {\sl Let $A$ be an abelian surface such that $NS(A)\simeq {\bf
Z}$,
linearly normally
embedded in ${\bf P}^8$ via a line bundle ${\cal L}$ giving a
polarization $\Theta$ in $NS(A)$ and let $\epsilon$ be any
non trivial point of order $2$ on $A$.
Then the map $\phi: X\to Y$ is an embedding.}\par

\vskip0.2truecm\noindent {\it Proof.} The polarization $\Theta$ is of type
$(d_1,d_2)$ with
$9=d_1\cdot d_2$, hence we have only the two cases $d_1=1, d_2=9$ and
$d_1=d_2=3$. In both cases $\epsilon\notin K(\Theta)$. In the
former case the assertion follows by theorem (3.3). Hence we consider only
the latter case, in
which $\Theta=3\theta$ in $NS(A)$, where $\theta$ is a principal
polarization.\par

By proposition (3.4), the morphism $\phi: X\to Y$ is birational.
Suppose two distinct points $z,
z'$ of $X$ are mapped to the same point $w$ by $\phi$.
By theorem (3.3, ii), $z$ and $z'$ do not
lie on $A$. Let
$F$ and
$F'$ be the two rulings of
$X$ through $z$ and $z'$ respectively, and let $x, x+\epsilon$ and
$y,y+\epsilon$ the pair of
points where $F$ and $F'$, respectively, meet $A$.
Then the points $w,x, x+\epsilon,y,y+\epsilon$
are coplanar in ${\bf P}^8$, hence we have a $4$-secant plane to $A$.
By proposition (3.2, iii),
there is an irreducible curve $C$, representing the polarisation
$\theta$, passing through $x,
x+\epsilon,y,y+\epsilon$. Consider the curve $C_\epsilon=t_\epsilon(C)$.
Since $\epsilon\not\in
K(\theta)=1$, it follows that $C_\epsilon\not=C$. On the other hand $x,
x+\epsilon,y,y+\epsilon$ belong to both $C$ and $C_\epsilon$.
Since $C\cdot C_\epsilon=C^2=2$, we
get a contradiction.\par

Suppose $d\phi$ is not injective at $z\in X$. Again $z\not\in A$. Let $F$
be the ruling of $X$
through $z$ and let $x, x'=x+\epsilon$ be the pair of
points where $F$ meets $A$. The same argument we made in the proof of
theorem (3.3, iii), shows
that there are two tangent lines $r$ and $r'$ to $A$ at $x$ and $x'$ which
lie in a
plane $\pi$ which is
therefore $4$-secant. Once more by proposition (3.2, iii) there is a curve $C$
representing
$\theta$, passing through
$x$ and
$x'$ and whose tangent cone at these points contains the
lines $r$ and $r'$. \par

Since $H^1(A,{\cal L}(-C))=
H^1(A,{\cal
O}_A(2C))=0$, the map $H^0(A,{\cal L})\to H^0(C,{\cal L}_{|C})$ is
surjective, i.e. $C$ is a curve of degree $6$ which spans a ${\bf P}^4$.
It follows that the linear system on $C$ cut out by the hyperplanes through
the plane  $\pi$ contains $2x+2x'$ and a residual $g^1_2$ which must be the
canonical $g^1_2$ of $C$, i.e.
${\cal L}_{|C}\simeq {\cal O}_C(2x+2x'+K_C)$.
As before $C\cap
C_\epsilon=\{x,x'\}$.
Hence $x+x'$, as a divisor on $C$, is linearly equivalent to $K_C+\eta$,
where $\eta$
is a suitable point of order two in $Pic^0(C)$. Therefore
$2x+2x'\equiv
2K_C$ and this yields that ${\cal L}_{|C}\simeq {\cal O}_C(3K_C)$. This
implies that
actually ${\cal L}\simeq {\cal O}_A(3C)$ (see [CFM, proposition (1.6)]). \par

This picks out exactly $81$ curves $C$ representing $\theta$ and
therefore $81$ rulings of $X$ on
which such a point $z$ can lie.
Note that these rulings form an orbit under the free action of the group
$K(\Theta)\simeq {\bf Z}_3^4$ on the set of rulings of $X$.
Hence the set of points $Z$ of $X$ where
$d\phi$ is not
injective
is finite. The image of $Z$ via $\phi$ is the set of singular points of
$Y$, which is
therefore also
finite. Moreover $Z$ is stable by the action of
$K(\Theta)\simeq {\bf Z}_3^4$ on
$X$. This implies that $Z$ is of order $81n$, so that $Z$ consists of
$81$ $n$-tuples of points,
each $n$-tuple lying on one of the aforementioned $81$ rulings of $X$. \par

Let us assume that $n\geq 1$. We consider the projection of $Y$ in ${\bf P}^6$
from a general line $R$ in
${\bf P}^8$. Since $Y$ has finitely many singularities the same holds for $Y'$.
The degree of
${\bf D}(X,6)$ is
$162$ and since each point of $Z$ clearly appears in ${\bf D}(X,6)$ with
multiplicity at least $2$,
we see that ${\bf D}(X,6)=2Z$. In particular $n=1$. As a consequence, we
also have that the
secant variety of
$Y$ cannot meet the general line $R$ of ${\bf P}^8$, hence it has
dimension $\nu\leq 6$, i.e. it has
dimension smaller than expected. This is excluded by a theorem of Scorza [S].
Hence we come
to a contradiction, which proves that $n=0$, thus proving that $\phi$ is an
embedding.\par

Scorza's argument in [S] is long and rather complicated. We give here, for
the reader's
convenience, a shorter version of it, adapted to our case.
Assume the secant variety of $Y$ has
dimension $\nu\leq 6$. By Terracini's lemma (see [LV, p. 18] or
[Z, prop. 1]), two
general
tangent spaces to $Y$ meet in
a subspace of dimension $6-\nu$. Actually we see that $\nu=6$.
Otherwise the general surface section $S$
is such that two general tangent spaces to it meet.
Then it is well known that
$S$ lies in a ${\bf P}^r$, with $r\leq 5$ (see [LV] or [CC]), a
contradiction.\par

Let us now make the projection $\psi$ of $Y$ to ${\bf P}^4$ from the
tangent space $\pi$ at a
general point
$y\in Y$. Since any other general tangent space to $Y$ meets $\pi$ at one
point,
we see that the
differential of $\psi$ has generic rank $2$. Hence $W=\psi(Y)$ is a
surface. \par

Notice that the general ruling of $Y$ does not meet $\pi$.
Otherwise the general
ruling would be contained in the span of
two general tangent spaces to $Y$, which is a
${\bf P}^6$. Then the whole $Y$ would be contained
in this ${\bf P}^6$, a
contradiction. Therefore,
since $Y$ is a
scroll, $W$ is also a scroll.
Notice that $W$ has only a $1$-dimensional
system of lines, otherwise it would be a plane, contrary to the fact that it
has to span a ${\bf
P}^4$. Let $R$ be a general line of $W$ and let $V$ be the closure of
$\psi^{-1}(R)$, which is a surface in a ${\bf P}^5$. Then there is
an irreducible
component
$V'$ of $V$ which is a scroll. The intersection of $V'$ with $A$ contains a
curve $C$
which is fixed by $t_\epsilon$. Hence $C$
represents $a\theta$
with $a$ positive and even. Furthermore $C$, as well as
$V'$, spans at most a
${\bf P}^5$. Consider the exact sequence:\par

$$0\to {\cal O}_A(-C)\otimes {\cal L}\to {\cal L}\to {\cal L}_{|C}\to 0$$

\noindent Since ${\cal O}_A(-C)\otimes {\cal L}$ represents
$(3-a){\theta}$, we see that
$H^1(A,{\cal O}_A(-C)\otimes {\cal L})=0$ (see [LB, p.66]).
Hence the restriction
map $H^0(A, {\cal
L})\to H^0(C,{\cal L}_{|C})$ is surjective. Then we must have:\par

$$6\geq h^0(C,{\cal L}_{|C})=h^0(A,{\cal L})-h^0(A,{\cal O}_A(-C)\otimes
{\cal L})=$$

$$=9-h^0(A,{\cal O}_A(-C)\otimes {\cal L})=9-h^0(A,{\cal O}_A((3-a)\theta))$$

\noindent i.e. $h^0(A,{\cal O}_A((3-a)\theta))\geq 3$.
Since $a$ is positive and even,
$h^0(A,{\cal O}_A((3-a)\theta))\leq h^0(A,{\cal O}_A(\theta))=1$, a
contradiction.$\diamondsuit$\par

\vskip 0.5truecm

\noindent {\centerline {\bf References.}}\par

\vskip 0.2truecm

\item{[BS]} M. Beltrametti, A. Sommese,  Zero-cycles and $k$-th
order embeddings of smooth projective surfaces, in:
Problems in the theory of surfaces and their
classification (F. Catanese, C. Ciliberto, M. Cornalba eds.), Symposia
Math. {\bf 32} (1991),
33-48.\par

\item{[CC]} L. Chiantini, C. Ciliberto, On a theorem of Palatini and
Terracini, in
preparation. \par

\item{[CFM]} C. Ciliberto, P. Francia, M. Mendes Lopes, Remarks on the
bicanonical
map for surfaces of general type, Math. Z., {\bf 224} (1997), 137-166.
\par

\item{[CH]} C. Ciliberto, K. Hulek, A remark on the geometry of elliptic
scrolls and
bielliptic surfaces, Manuscripta Math. {\bf 95} (1998), 213-224.\par

\item{[EH} D. Eisenbud, J. Harris, Curves in projective space, Montreal
Univ. Press, 1982.\par

\item{[F]}W. Fulton, Intersection theory, Ergebnisse der Mathematik,
Springer Verlag, 1984.\par

\item{[LB]} H. Lange, Ch. Birkenhake, Complex abelian varieties,
Grundlehren der Mathematik,
Springer Verlag, 1992.\par

\item{[LV]} R. Lazarsfeld, A. van de Ven, Topics in the geometry of
projective space, DMV,
Birk\"auser, 1984.\par

\item{[R]} E. Rogora, Varieties with many lines, Manuscripta Math., {\bf 82}
(1994), 207-226.\par

\item{[S]} G. Scorza, Determinazione delle variet\`a a tre dimensioni di
$S_r$ $(r\geq 7)$
i cui
$S_3$ tangenti si tagliano a due a due,
Rend. Circ. Mat. Palermo, {\bf 25} (1908), 193-204.\par

\item{[Z]} F. Zak, Projections of algebraic varieties, Math. USSR Sbornik,
{\bf 44}
(1983), 535-544.\par

\vskip 4truecm

\noindent Authors' addresses:

\vskip 0.5truecm

\noindent Ciro Ciliberto\par
\noindent Universit\`a di Roma Tor Vergata\par
\noindent Dipartimento di Matematica\par
\noindent Via della Ricerca Scientifica\par
\noindent I 00173 Roma (Italy)\par
\noindent e-mail: Cilibert@axp.mat.uniroma2.it\par

\vskip 0.5truecm

\noindent Klaus Hulek\par
\noindent Institut f\"ur Mathematik\par
\noindent Universit\"at Hannover\par
\noindent Post fach 30060\par
\noindent D 3167 Hannnover (Germany)\par
\noindent e-mail: Hulek@math.uni-hannover.de\par

\end